\documentclass{birkjour}

\usepackage{amsmath,amssymb,amsthm,color}
\usepackage[utf8]{inputenc}
\usepackage{amsopn}
\usepackage{latexsym}
\usepackage[all,cmtip]{xy}
\usepackage{tikz}
\usepackage{float}

\theoremstyle{plain}
 \newtheorem{teo}{Theorem}[section]

 \newtheorem*{defi}{Definition} 
 
\theoremstyle{definition}
\newtheorem{ex}[teo]{Example} 
\newtheorem*{remark}{Remark}

\newcommand{\nc}{\newcommand}

 \nc{\ad}{\operatorname{ad}} 
 \nc{\Ad}{\operatorname{Ad}}
 \nc{\coad}{\operatorname{coad}} 
 \nc{\ct}{\operatorname{T}}
 \nc{\rank}{\operatorname{rank}} 
 \nc{\Irr}{\operatorname{Irr}}
 \nc{\End}{\operatorname{End}} 
 \nc{\Aut}{\operatorname{Aut}}
 \nc{\Inn}{\operatorname{Inn}} 
 \nc{\Der}{\operatorname{Der}}
 \nc{\Dera}{\operatorname{Dera}} 
 \nc{\Auto}{\operatorname{Auto}}
 \nc{\GL}{\operatorname{GL}}
 \nc{\SL}{\operatorname{SL}}
 \nc{\coord}{\operatorname{coord}}
 \renewcommand{\span}{\operatorname{span}}
 \nc{\codim}{\operatorname{codim}}
  \nc{\pr}{\operatorname{pr}}

 \newcommand{\R}{\mathbb R}

\newcommand{\mn}{\mathfrak n }

\newcommand{\mz}{\mathfrak z }
\newcommand{\mk}{\mathfrak k }
\newcommand{\mc}{\mathfrak c }

\newcommand{\mh}{\mathfrak h }
\newcommand{\ma}{\mathfrak a }
\newcommand{\mgg}{\mathfrak g }

\renewcommand{\sl}{\mathfrak{sl} }

\newcommand{\talpha}{ \tilde{\alpha} }

\nc{\rmc}{\textrm{\rm C}}
\nc{\rmd}{\textrm{\rm D}}
\nc{\rad}{\textrm{\rm Rad}}

\newcommand{\lra}{\longrightarrow}

\begin{document}
\title[On generalized $G_2$-structures and $T$-duality]{On generalized $G_2$-structures and $T$-duality}

\author{Viviana del Barco}
\email{delbarc@ime.unicamp.br}

\address{UNR-CONICET (Argentina) and IMECC-UNICAMP (Brazil)}

\author{Lino Grama}
\email{linograma@gmail.com}
\address{IMECC-UNICAMP}

 \date{\today}

\begin{abstract} This is a short note on generalized $G_2$-structures obtained as a consequence of a $T$-dual construction given in \cite{dBGS}. Given classical $G_2$-structure on certain seven dimensional manifolds, either closed or  co-closed, we obtain integrable generalized $G_2$-structures which are no longer an usual one, and with non-zero three form in general. In particular we obtain manifolds admitting closed generalized $G_2$-structures not admitting closed (usual) $G_2$-structures. 
\end{abstract}

\thanks{V. del Barco supported by FAPESP grant 2015/23896-5. \\
L. Grama supported by FAPESP grant 2016/22755-1 and 2012/18780-0.}



\maketitle

\section{Introduction}

Generalized $G_2$-structures were introduced by Witt \cite{Wit} with the intention of generalizing the classical concept of $G_2$-structure defined by a 3-form $\varphi$ \cite{Joy00}. A generalized $G_2$-structure on a 7-dimensional  differentiable manifold  $M$ is a reduction from the structure group $\R^*\times Spin(7,7,)$ of the principal bundle $TM\oplus T^*M$  to $G_2\times G_2$. This reduction determines a generalized metric and a 2-form on $M$, which split $TM\oplus T^*M$ into submodules with positive and negative definite metrics. Therefore there is associated a pair of spinors $\Psi_\pm$ in the irreducible spin representation $\Delta=\R^8$ of $Spin(7)$. The  $G_2\times G_2$ invariant tensor $\Psi_+\otimes \Psi_-\in \Delta\otimes \Delta$ can be considered  as a differential  form on $M$, so it induces elements $(\Psi_+\otimes \Psi_-)^{even/odd}$ corresponding to the even and odd degrees \cite{Wit} . Up to a $B$-field transformation and a dilaton we have (see also \cite{HuHu}) 
\begin{eqnarray*}
\rho &=& \Psi_+\otimes \Psi_-^{even} =  c-c\star\varphi+s \star( \alpha\wedge\star\varphi)-s \alpha\wedge \varphi-s\star\alpha,\\
\hat\rho &=& \Psi_+\otimes \Psi_-^{odd} =s\alpha-c\varphi-s \star( \alpha\wedge\varphi)-s \alpha\wedge \star \varphi+c\frac17 \varphi\wedge\star\varphi,
\end{eqnarray*} where $\alpha$ is a unit 1-form,  $\varphi$ is a 3-form and the parameters $s,c$ correspond, respectively, to the sine and cosine of the angle between the spinors $\Psi_{\pm}$; in particular  $s^2+c^2=1$.

Let $H$ be a 3-form on $M$.  A generalized $G_2$-structure defined by the spinors $\rho$, $\hat \rho$ as above, is called {\em strongly integrable} with respect to $H$ if 
\begin{equation}
d_H\rho=d_H\hat\rho=0,\label{eq.strong}\end{equation}
where $d_H\cdot=d\cdot+H\wedge$ is the twisted operator of $d$. The generalized $G_2$-structure is called {\em weakly integrable of odd (resp. even) type} if 
\begin{equation}
\label{eq.weakly}
d_H\hat\rho=\lambda\rho,\quad \mbox{(resp. }d_H\rho=\lambda\hat \rho\mbox{)}.\end{equation}
for some non-zero constant $\lambda$. The real number $\lambda$ is called the Killing number.\smallskip

An usual $G_2$-structure $\varphi$ on a 7-manifold induces a generalized $G_2$-structure on $M$. In this case $\Psi_+=\Psi_-$ and $s=0$, therefore the even and odd spinors are given by
\begin{equation}\label{usualg2}
\rho = 1-\star\varphi,\quad \hat\rho = -\varphi+dV.
\end{equation} When $H$ is zero, we see that the generalized structure is strongly integrable \eqref{eq.strong} if and only if the usual $G_2$ is parallel, that is, $\varphi$ is closed and co-closed. Instead,  weak integrability of odd type cannot occur for usual $G_2$-structures, independently of $H$ (see \cite[Page 288]{Wit}). If we let $H$ to be any closed 3-form, the only compact strongly integrable generalized $G_2$ are the usual parallel $G_2$-structures \cite{Wit}. Fino and Tomassini  gave the first example of a compact strongly integrable generalized $G_2$-structure for $H$ non-closed \cite{FT}.

Making an analogy with the classical case we define
a generalized $G_2$-structure with spinors $\rho$ and $\hat\rho$ to be a closed (resp. co-closed) structure if 
\begin{equation}
\label{eq.closed}
d_H\hat\rho=0,\quad \mbox{(resp. }d_H\rho=0\mbox{)}.\end{equation}
Given a closed $G_2$-structure, the generalized structure associated to $\varphi$ is closed with respect to any 3-form $H$ such that $H\wedge \varphi=0$. To the contrary, if $\varphi$ is co-closed, then the generalized structure is co-closed only for $H=0$. In fact by \eqref{usualg2}, $d_H\rho=H-H\wedge \star\varphi-d\star\varphi=0$ if and only if $H=0$.

Witt himself was interested in the relation of generalized $G_2$-structures defined on $T$-dual manifolds. Recall that $T$-duality is a concept introduced by Bouwknegt, Evslin, Hannabuss and Mathai \cite{BEM}, \cite{BHM} for manifolds having the structures of principal torus bundles. According to \cite{BEM}, if $H$ and $H^\vee$ are closed 3-forms on $T$-dual manifolds $M$ and $M^\vee$ are $T$-dual then there exists an isomorphism $\tau$ between the differential complexes $(\Omega_{T^k}^{\bullet} (M),d_H)$ and $(\Omega_{T^k}^{\bullet} (M^\vee),d_{H^\vee})$. Such isomorphism $\tau$ satisfies
$$ d_{H^\vee}\tau (\rho)=\tau (d_H \rho). $$
Therefore the integrability conditions are preserved by $T$-duality, thus generalized $G_2$-structures on $M$ integrable with respect to $H$, induce generalized $G_2$-structures on $M^\vee$, integrable with respect to $H^\vee$. And vice-versa.
\smallskip

The aim of this short note is to contribute with further examples of integrable generalized $G_2$-structures following the construction of $T$-duals the authors developed in \cite{dBGS}, together with Leonardo Soriani.

We start with an usual $G_2$-structure on a seven dimensional manifold seen as a generalized $G_2$ with  respect  to certain three forms $H$. In the dual manifold to the given one we obtain an integrable generalized $G_2$ which is no longer an usual one, and with non-zero three form in general. Our  previous work focuses on invariant structures on compact quotients of solvable Lie groups by discrete subgroups. So this framework maintains here and we develop the  examples at the Lie algebra level. 

In \cite{dBGS} we indicated how duality would contribute with the study of symplectic structures by dualizing generalized complex structures. A similar spirit is pursued here, but in this case we work directly with the spinors defining the generalized structures. In this manner we are able to present manifolds admitting closed generalized $G_2$-structures not admitting closed (usual) $G_2$-structures. This note was motivated by the papers of Witt, Fino and Tomassini \cite{FT,Wit}.

\section{Lie algebras and infinitesimal duality} 

We briefly recall the concept of duality of Lie algebras \cite{dBGS}.

Let $\mgg$ be a Lie algebra together with a closed 3-form $H$. Let $\ma$ be an abelian ideal of  $\mgg$, we say that the triple $(\mgg,\ma,H)$ is admissible if $H(x,y,\cdot)=0$ for all $x,y\in\ma$. Notice that when $\dim \ma=1$ then any closed 3-form gives an admissible triple.

 In this case, denote $\mn$ the quotient Lie algebra and $q:\mgg\lra \mn$  and $q^\vee:\mgg^\vee\lra \mn$ the quotient maps. The subspace $\mc$ of $\mgg\oplus\mgg^\vee$
$$\mc=\{(x,y)\in \mgg\oplus\mgg^\vee: q(x)=q^\vee(y)\} $$ is a Lie subalgebra and the  following diagram is  commutative
$$\xymatrix{
    &\mc \ar[ld]_{p}  \ar[rd]^{p^\vee}& \\
\mgg\ar[dr]_{q}&  & \mgg^\vee\ar[ld]^{q^\vee}.\\
               & \mn  & }$$
Here $p$ and $p^\vee$ are the projections over the first and second component, respectively. A 2-form $F\in \Lambda^2\mc^*$ is said to be {\em non-degenerate in the fibers} if for all $x\in \mk=\{(x,0)\in\mc:x\in \ma\}$, there exists some $y\in \mk^\vee=\{(0,y)\in\mc:y\in \ma^\vee\}$ such that $F(x,y)\neq 0$. Such an $F$ exists if and only if $\dim \ma=\dim\ma^\vee$. 

\begin{defi} Two  admissible triples $(\mgg,\ma,H)$ and $(\mgg^\vee,\ma^\vee,H^\vee)$ are said to be {\em dual} if $\mgg/ \ma \simeq\mn\simeq \mgg^\vee/\ma^\vee$ and there exist a 2-form $F$ in $\mc$ which is non-degenerate in the fibers such that $p^*H-p^{\vee*}H^\vee=dF$.
\end{defi}

The dual of a given admissible triple always exists and its construction is described in the following result.

\begin{teo}\cite{dBGS}\label{teo.infdual} Let $(\mgg,\ma,H)$ be an admissible triple with $\ma$ a central ideal and let $\{x_1,\ldots,x_m\}$ be a basis  of  $\ma$. Define 
\begin{itemize}
\item $\Psi^\vee=(\iota_{x_1}H,\ldots,\iota_{x_m}H)$, 
\item $\mgg^\vee=(\mgg/\ma)_{\Psi^\vee}$ and
\item $H^\vee=\sum_{k=1}^m z^k\wedge dx^k+\delta$ where $\{z_1,\ldots,z_m\}$ is a basis of $\ma^\vee$ and $\delta$ is the basic component of $H$.
\end{itemize}
Then $(\mgg^\vee,\ma^\vee,H^\vee)$ is an admissible triple and is dual to $(\mgg,\ma,H)$. 

Conversely, if $(\mgg^\vee,\ma^\vee,H^\vee)$ is dual to $(\mgg,\ma,H)$, then there exist a basis $\{x_1,\ldots,x_m\}$ of  $\ma$ and a basis $\{z_1,\ldots,z_m\}$ of $\ma^\vee$ such that the formulas above hold.
\end{teo}

Here $(\mgg/\ma)_{\Psi^\vee}$ denotes the central extension of $\mgg/\ma$ by the closed 2-form $\Psi^\vee$ (see \cite{dBGS} for details).

\bigskip

$G_2$ and $SU(3)$-structures on Lie algebras were treated by several authors. In our context, the main references are the work of Conti, Fernandez, Fino, Manero, Raffero and Tomassini (see \cite{CoFe,FT,Man} and references therein).
A $G_2$-structure on a Lie algebra $\mgg$ is a non-degenerate 3-form $\varphi$ such that is some adapted basis $\{e^1,\ldots, e^7\}$ of the dual $\mgg^*$  of $\mgg$, it is written as
 $$\varphi=e^{127}+e^{347}+e^{567}+e^{135}-e^{236}-e^{146}-e^{245}.$$
 
To continue, we consider Lie algebras $\mgg$ endowed with $G_2$-structures and with non-trivial center. We consider central ideals in $\mgg$ and their dimensions is what we call the dimension of the fiber, having in mind a possible torus bundle structure on the Lie group associated to $\mgg$. After fixing a closed 3-form giving integrability of the usual $G_2$-structure, we compute the dual triple and the generalized $G_2$-structure arising on it.

\subsection{One dimensional fiber}

Let $\mgg$ be a Lie algebra with non-trivial center and let $\ma$ be a one dimensional central ideal. Fix  a generator $x\in \ma$ and let $\mgg=\ma\oplus \mh$ be an orthogonal decomposition with respect to the metric induced by $\varphi$. Let  $\varphi$ be a $G_2$-structure on $\mgg$, then if $\alpha$ is the dual element to $x$ we have
$$\varphi=\alpha\wedge \omega+\psi_+,\mbox{ with }\iota_x\omega=0\mbox{ and }\iota_x\psi_+=0.$$
Up to a normalization of coefficients, the forms $(\omega,\psi)$ define an $SU(3)$-structure on $\mn=\mgg/\ma$ \cite{ApSa,CoFe}.

The spinors associated to the generalized $G_2$-structure induced by $\varphi$ are given by \eqref{usualg2}:
\begin{equation}\label{dualspinorsfiber1}
\rho = 1-\frac12 \omega^2-\psi_-\wedge \alpha,\quad \hat\rho = -\alpha\wedge\omega-\psi_++dV.
\end{equation}

Any closed 3-form $H$ in $\mgg$ makes $(\mgg,\ma,H)$ an admissible triple, since $\dim\ma=1$. According to Theorem \ref{teo.infdual}, the dual triple is $(\mgg^\vee,\ma^\vee,H^\vee)$ where
$\mgg^\vee$ is the central extension of $\mn$ by $\iota_xH$. Explicitly, $\mgg^\vee=\R z\oplus \mn$ and the Lie bracket of $\mgg^\vee$ satisfies
$$[u,v]=[u,v]_{\mn}+(\iota_xH)(u,v) \,z \;\mbox{ and }\;[z,u]=0, \mbox{  for all }u,v\in\mn.$$
Let $\talpha$ be the 1-form such that $\talpha(z)=1$ and $\talpha(\mn)=0$, then the closed 3-form on $\mgg^\vee$ is $H^\vee=\talpha\wedge d\alpha +\delta$ where $\delta$ the basic part of $H$. Note that $H\neq 0$ if and only if $d\alpha\neq 0$ which means that $\R x$ is not a direct factor of $\mgg$. The 2-form $F\in \Lambda^2\mc^*$ giving the duality on the correspondence  space is given by $F=\tilde \alpha\wedge\alpha$.

The dual spinors $\rho^\vee$ and $\hat\rho^\vee$ are given by (see \cite{CG})
$$\rho^\vee=\iota_xe^F\rho,\qquad \hat\rho^\vee=\iota_xe^F\hat\rho. $$ Explicitly, we obtain
\begin{equation}\label{dualsp}
\rho^\vee= -\tilde\alpha+\psi_-+\frac12\tilde\alpha\wedge\omega^2,\quad \hat\rho^\vee =-\omega+\tilde\alpha\wedge\psi_++\frac16\omega^3.
\end{equation}

\begin{remark}
The spinors above correspond to a generalized $G_2$-structure associated to the $SU(3)$-structure $(\omega,\psi_+)$, when one imposes the angle of the spinors to be $\pi/2$. These structures where considered in \cite{FT}.
\end{remark}

\begin{remark}
Notice that the dual of an usual $G_2$-structure is never an usual $G_2$-structure, but a pure generalized $G_2$.
\end{remark}

\begin{ex} \label{ex1} Let $\mgg$ be the Lie algebra spanned by $\{e_1, \ldots, e_7\}$ and satisfying the bracket relations
$$[e_1,e_7]=-e_3,\;[e_1,e_5]=-e_4,\;[e_2,e_7]=-e_4,\;[e_1,e_3]=-e_6,$$
and zero in the other cases. The Lie differential in the dual basis is $de^3=e^{17}$, $de^4=e^{15}+e^{27}$ and $de^6=e^{13}$. Consider the central ideal $\ma=\R  e_6$.

The 3-form 
$\varphi=e^{127}+e^{347}+e^{567}+e^{135}-e^{236}-e^{146}-e^{245}$ is a closed $G_2$-structure \cite[Example 1.5]{Man}. Notice that $\varphi=e^6\wedge\omega+\psi_+$ with $\omega=-e^{14}-e^{23}-e^{57}$ and $\psi_+=e^{127}+e^{347}+e^{135}-e^{245}$. The pair $(\omega,\psi_+)$ define an $SU(3)$-structure on $\mgg/\ma$

 We consider all possible closed 3-forms $H$ such that $d_H\hat\rho=0$. Canonical computations give that $H$ is of the form
\begin{eqnarray}
H&=&a_1(e^{134}+e^{267}+e^{123}+e^{357})
+a_3(e^{136}+e^{137}+e^{145}-e^{247})\nonumber\\
&&\;\;+a_2(e^{126}-e^{237})+a_4(e^{156}-e^{357})+a_5(e^{157}+e^{134}+e^{267}+e^{357})\nonumber\\
&&\;\;+a_6(e^{167}-e^{257})+a_7(e^{367}-e^{457})+a_8(e^{146}+e^{236})\label{Hcompatible}\\
&&\;\;+a_9(e^{347}+e^{567})+a_{10}(e^{124}+e^{257})+a_{11}(e^{125}+e^{137})\nonumber\\
&&\;\;+a_{12}e^{127}+a_{13}e^{135}+a_{14}
e^{245}+a_{15}(e^{145}+e^{235}),\quad a_i\in \R.\nonumber
\end{eqnarray}

For any $H$ as in \eqref{Hcompatible}, the triple $(\mgg,\ma,H)$ is a compatible triple and $d_H\hat\rho=0$. Thus $\rho,\hat\rho$ define a closed generalized $G_2$-structure with respect to $H$. Now we describe the dual triple.

The Lie algebra $\mgg^\vee$ has a dual basis $\{f^1, \ldots,f^7\}$ such that the Lie algebra differential is $$
\begin{array}{rcl}
df^3&=&f^{17},\\
df^4&=&f^{15}+f^{27},\\
df^6&=&\iota_{e_6H}=-a_6 f^{17}-(a_1+a_5) f^{27}-a_7 f^{37}-a_9 f^{57}\\
&&\;\;+a_2 f^{12}+a_3  f^{13}+a_4 f^{15}+a_8 (f^{14}+f^{23}).\end{array}
$$ Here we identify $e_i$ with $f_i$ for $i=1, \ldots, 5, 7$, taking into account that $\mgg/\ma\simeq \mgg^\vee/\ma^\vee$. The dual 3-form is 
\begin{eqnarray*}
H^\vee &=&f^{136}+a_1(f^{134}+f^{123}+f^{357})-a_2 f^{237}+a_3(f^{137}-f^{247}+f^{145})\\
&&\;\;-a_4 f^{357}+a_5(f^{157}+f^{134}+f^{357})-a_6 f^{257}-a_7 f^{457}+a_9 f^{347}\\
&&\;\;+a_{10}(f^{124}+f^{257})+a_{11}(f^{125}+f^{137})+a_{15}(f^{145}+f^{235})\\
&&\;\;+a_{12}f^{127}+a_{13}f^{135}+a_{14}
f^{245}.
\end{eqnarray*}
and the dual spinor $\hat\rho^\vee$ is, by \eqref{dualsp}, 
$$\hat\rho^\vee=-\omega+f^6\wedge \psi_+=f^{14}+f^{23}+f^{57}+f^{2456}+f^{1267}+f^{3467}-f^{1356}.$$ One can check directly that $d_{H^\vee}\hat\rho^\vee=0.$

Depending on the coefficients $a_i$ of $H$ in \eqref{Hcompatible}, we reach non-isomorphic Lie algebras. For instance, if $H=e^{136}+e^{137}-e^{247}+e^{145}$ ($a_3=1$ while the others are zero), then $\mgg^\vee\simeq \mgg$. Meanwhile $\mgg^\vee$
%=(0,0,17,15+27,0,14+23,0)$
has one dimensional center for $H=e^{146}+e^{236}$. If $H=0$ then the only non-trivial Lie brackets of $\mgg^\vee$ are $[f_1,f_7]=-f_3$ and $[f_1,f_5]=-f_4=[f_2,f_7]$.

When $H$ is as in \eqref{Hcompatible} and satisfies $a_4\cdot a_7=0$ and $a_3^2+a_8^2>0$  then $(\iota_{e_6}\sigma)^3=0$ for any closed 3-form in the dual Lie algebra $\mgg^\vee$, so these Lie algebras do not admit closed $G_2$-structures, but they admit closed generalized $G_2$, as we showed above.

\end{ex}

\begin{ex} The Lie algebra $\mgg$ of dimension 7 with non-zero Lie brackets $[e_2,e_5]=-e_6$, $[e_4,e_5]=e_7$ admits closed $G_2$-structures \cite{CoFe}, such as, for instance, $\varphi=e^{127}+e^{347}+e^{567}+e^{135}-e^{146}-e^{236}-e^{245}$. Consider the central ideal $\ma$ spanned by $e_7$
and let $H$ be the 3-form
\begin{eqnarray}
H&=&a_1(e^{124}-e^{456})+a_2(e^{125}-e^{345})-a_3(e^{134}-e^{156})+a_4\,e^{135}+\nonumber\\
&&\quad a_5(e^{145}-e^{235})+ a_6(e^{145}+e^{246})+a_7(e^{234}-e^{256})+a_8\,  e^{245},\label{HTF} \end{eqnarray}
where $a_i$ are real coefficients. Thus $H$ is closed and satisfies $d_H\hat\rho=d_H(-\varphi+dV)=-d\varphi-H\wedge \varphi=0$. Therefore, we have a closed generalized $G_2$-structure with respect to $H$. 

Since $\iota_{e_7}H=0$, the dual Lie algebra $\mgg^\vee$ is defined by the only non-zero bracket relation $[f_2,f_5]=-f_6$ (as before we identify $e_i$ with $f_i$ for $i=1, \ldots, 6$). The dual 3-form is \begin{eqnarray*}
H^\vee&=&-f^{457}+a_1(f^{124}-f^{456
})+a_2(f^{125}-f^{345})-a_3(f^{134}-f^{156})+\\
&&\;\; a_4\,f^{135}+a_5(f^{145}-f^{235})+a_6(f^{145}+f^{246})+a_7(f^{234}-f^{256})\\
&&\;\;+a_8\,  f^{245}. 
\end{eqnarray*}
The duality preserving the integrability, implies that  $d_{H^\vee}\hat\rho^\vee=0$. Thus we re-obtained Example 5.1. of Fino and Tomassini in \cite{FT}.

\end{ex}

\subsection{Fiber of dimension $2$} When $\ma$ has dimension greater than one, there is no global expressions as \eqref{dualspinorsfiber1} for the dual spinors, so their computations need to be done by hand in each case. 

\begin{ex} Let $\mgg$ be the Lie algebra such that the Lie algebra differential on a dual basis $\{e^1,\ldots,e^7\}$ is
$$de^1=e^{35}+e^{46},\;de^3=e^{67},\;de^4=e^{57},\;de^5=e^{47},\;de^6=e^{37},\;de^2=de^7=0. $$ This Lie algebra is solvable  and unimodular, and it admits a co-closed $G_2$ form. Indeed, 
  $$\varphi=e^{127}+e^{347}+e^{567}-e^{136}-e^{145}-e^{235}+e^{246} $$ satisfies
  $$\star \varphi= e^{1234}+e^{1256}+e^{3456}+e^{1357}-e^{1467}-e^{2367}-e^{2457}  $$
  and   $d\star \varphi=0$ (see \cite[Example 2.6]{Man}).
  
  Denote  $\mh$ the Lie algebra with basis $\{e_1, \ldots,e_6\}$ and unique non-zero differential $de^1=e^{35}+e^{46}$. Then $\mgg$ is the  one dimensional extension of $\mh$ by the derivation $D$ defined as $De_{3+i}=e_{6-i}$, $i=0,\ldots,3$.
  
The spinors associated to this $G_2$-structure are given in \eqref{usualg2}:
\begin{eqnarray*}
\rho &=&1- (e^{1234}+e^{1256}+e^{3456}+e^{1357}-e^{1467}-e^{2367}-e^{2457})  \\
\hat\rho &=&-(e^{127}+e^{347}+e^{567}-e^{136}-e^{145}-e^{235}+e^{246})+e^{1234567}.
\end{eqnarray*}

The triple $(\mgg,\ma,H=0)$ with  $\ma=\mz(\mgg)=\span\{e_1,e_2\}$ is a compatible triple. And the generalized $G_2$-structure is co-closed with respect to $H=0$ (recall that co-closed $G_2$-structures are co-closed generalized $G_2$-structures only for $H=0$). We shall compute the dual triple and structure.

Since $\iota_{e_2}H=\iota_{e_1}H=0$, the dual Lie algebra $\mgg^\vee$ is the trivial extension of $\mn=\mgg/\ma$, thus it has a dual basis $\{f^1,\ldots ,f^7\}$ such that the differentials account to
$$df^3=f^{67},\;df^4=f^{57},\;df^5=f^{47},\;df^6=f^{37},\;df^1=df^2=df^7=0. $$ 
One has that $\mgg^\vee$ is isomorphic to $\R^2\oplus (\R\ltimes_D \R^4)$, with $D$  as above, induced to $\R^4$.
The dual 3-form is $H^\vee=f^{135}+f^{146}$.

The dual spinor $\rho^\vee$ is 
\begin{eqnarray*}
\rho^\vee&=&-(f^{34}+f^{56}+f^{12})+f^1(f^{367}+f^{457})+
f^2(f^{357}-f^{467})+f^{123456} \\
&=&-(f^{12}+f^{34}+f^{56})+(f^{136}+f^{145}+f^{235}
-f^{246})f^7+f^{123456}
\end{eqnarray*}
One can verify that $d_H\rho^\vee=0$. 

The Lie algebra $\mgg^\vee$ also admits a co-closed usual $G_2$-structure since it is of the form $\R f_7 \ltimes_D \R^6$ for the derivation given above and $\omega=f^{12}+f^{34}+f^{56}$, $\psi_+=f^{135}-f^{146}-f^{236}-f^{246}$ define a half flat structure on $\R^6$. So  Proposition 2.1. in \cite{Man} give co-closed $G_2$-structures on $\mgg^\vee$.
\end{ex}

\bibliographystyle{plain}
\bibliography{biblio}

\end{document}